\input amstex
\documentstyle{amsppt}
\NoBlackBoxes
\input bull-ppt
\keyedby{bull238/jxs} 

\topmatter
\cvol{26}
\cvolyear{1992}
\cmonth{Jan}
\cyear{1992}
\cvolno{1}
\cpgs{1-2}
\title Editors' remarks \endtitle
\author Morris W. Hirsch and Richard S. Palais\endauthor
\shortauthor{M. W. Hirsch and R. S. Palais}
\address Department of Mathematics, University of 
California at Berkeley,
Berkeley, California 94720\endaddress
\address Department of Mathematics, Brandeis University, 
Waltham,
Massachusetts 02154-9110\endaddress
\endtopmatter

\document

The two papers that follow are controversial---in two 
senses.  First,
the authors express opposing and strongly worded views on 
what
``complexity theory'' should be.  Second, the decision to 
open
the Research-Expository Papers section of the {\it 
Bulletin\/}
as a forum for such a debate may also be considered 
controversial; and,
in fact, the wisdom of this decision was the subject of 
some dispute
within the editorial board.  It was made by us jointly as 
present
and past chairs of the editorial board.

As to the controversy in the first sense, we let the 
papers speak
for themselves.  However, since their publication is a 
precedent
of sorts, we feel it is important to clarify our general 
attitudes
toward articles of a controversial nature.

As mathematicians we have the good fortune to be able to 
settle
in a straightforward and objective way one sort of 
controversy 
which, in other disciplines, often leads to quite 
rancorous 
disputes.  While there are occasional disagreements over 
the correctness
of a paper, the strictly logical nature of mathematical 
proof usually
permits a quick resolution of such issues that is agreed 
to by all 
sides.  But this should not blind us to other mathematical 
controversies that are less objective in their nature, 
and not
so easily settled.

For example, we have probably all heard the story that 
some 
mathematicians felt it was scandalous for Cantor to claim 
that,
in demonstrating that the algebraic numbers were countable
while the real numbers were not, he had given a new proof 
of the
existence of transcendental numbers.  After all, his 
proof gave
no way to construct even a single transcendental number. 
 Echoes
of this controversy are heard down to the present day in 
the now
somewhat muffled debate over ``Constructivism versus 
Classical 
Mathematics.''  Similarly, we read that Hilbert's 
approach to 
Invariant Theory, using his Basis Theorem and other 
nonconstructive,
abstract methods, provoked controversy in a mathematical 
world still
steeped in the concrete methods of the classical 
tradition, where 
solving a particular problem in Invariant Theory had 
always meant
exhibiting a specific basis for the invariants.  Other 
controversies
include debates over the status of infinitesimals, 
irrationals,
imaginary numbers, large cardinals; the proper treatment 
of geometry,
logic, set theory, foundations of mathematics; the role 
of computer
science; the use of mathematics in the social sciences; 
and perennial
issues in mathematical education.

And not all such controversies are ancient history! 
 Fifteen years ago
there was a sharp controversy over purported excesses in 
the applications
of Catastrophe Theory, and currently there is a similar 
controversy 
concerning what some see as an overselling and 
overpopularization of 
``fractals'' and ``chaos.''  Another simmering debate has 
grown out of
the current renewal of the on-again, off-again love 
affair between
Mathematics and Theoretical Physics.  We have learned to 
accept that
different standards of mathematical rigor may be 
appropriate when
mathematics is being used as a tool to gain new insights 
about the 
physical world.  But what standards should we apply to 
judge a paper
that uses nonrigorous or semirigorous methods from 
physics to 
suggest important new insights into our own mathematical 
world, 
particularly if those insights seem beyond the reach of 
current rigorous
mathematics?

Especially because such questions cannot always be 
answered by logical
principles alone, we believe that it is important for 
mathematicians to
confront them.  Even when rational discussion and debate 
does not 
completely resolve differences, at least it may clarify 
the issues.

Traditionally debate about issues of this sort has been 
carried on
in nonscholarly journals, and for questions that are less 
weighty or more
transitory in significance this is appropriate.  We 
certainly have no
intention to open these pages to emotional debate over 
whether the C
programming language is better or worse than Pascal!  But 
when a controversial
matter comes up that is of serious concern and long-term 
significance to
the mathematical community, and so deserving of careful 
debate, then such
a debate belongs in an archival journal.  This does {\it 
not\/} mean 
we are inviting authors to submit some new category of 
``controversial
issue'' paper to the Research-Expository Papers section. 
 On the 
contrary, as always, any paper will judged on its 
intrinsic interest and
merits, and controversial papers will no doubt have to 
jump through a 
few extra hoops.  What we are saying is that we will not 
reject a paper
solely because the ideas presented in it may not be 
universally accepted
or subject to mathematical proof or disproof.

\enddocument